
\tolerance=10000
\magnification=1200
\raggedbottom

\baselineskip=15pt
\parskip=1\jot

\def\sk{\vskip 3\jot}

\def\heading#1{\vskip3\jot{\noindent\bf #1}}
\def\label#1{{\noindent\it #1}}
\def\QED{\hbox{\rlap{$\sqcap$}$\sqcup$}}


\def\ref#1;#2;#3;#4;#5.{\item{[#1]} #2,#3,{\it #4},#5.}
\def\refinbook#1;#2;#3;#4;#5;#6.{\item{[#1]} #2, #3, #4, {\it #5},#6.} 
\def\refbook#1;#2;#3;#4.{\item{[#1]} #2,{\it #3},#4.}


\def\({\bigl(}
\def\){\bigr)}


\def\al{\alpha}

\def\de{\delta}
\def\ep{\varepsilon}
\def\ze{\zeta}

\def\th{\vartheta}

\def\ta{\tau}

\def\ph{\phi}

\def\Th{\Theta}

\def\Me{\Omega}


\def\calM{{\cal M}}

\def\calU{{\cal U}}
\def\calV{{\cal V}}
\def\calW{{\cal W}}
\def\calX{{\cal X}}
\def\calY{{\cal Y}}


\def\XY{X\cdot Y}
\def\Eq{H(X,Y\mid\XY)}
\def\xy{x\cdot y}

\def\gcd{{\rm gcd}}

{
\pageno=0
\nopagenumbers
\rightline{\tt multiplication.tex}
\vskip1in

\centerline{\bf The Average Amount of Information}
\centerline{\bf Lost in Multiplication}
\vskip0.5in

{\baselineskip=12pt
\centerline{Nicholas Pippenger}
\centerline{\tt njp@princeton.edu}
\vskip0.25in

\centerline{Department of Computer Science} 
\centerline{Princeton University}
\centerline{35 Olden Street}
\centerline{Princeton, NJ 08540 USA}
\vskip1in
}

\noindent{\bf Abstract:}
We show that if $X$ and $Y$ are integers independently and uniformly distributed in
the set $\{1, \ldots, N\}$, then the information lost in forming their product
(which is given by the equivocation
$\Eq$), is $\Th(\log\log N)$.
We also prove two extremal results regarding cases in which $X$ and $Y$ are 
not necessarily independently
or uniformly distributed.
First, we note that
the information lost in multiplication can of course be $0$.
We show that the  condition
$\Eq = 0$ implies $2\log_2 N - H(X, Y) = \Me(\log\log N)$.
Furthermore, if $X$ and $Y$ are independent and uniformly distributed on disjoint sets of primes,
it is possible to have $\Eq = 0$ with $\log_2 N - H(X)$ and $\log_2 N - H(Y)$
each $O(\log\log N)$.
Second, we show that however $X$ and $Y$ are distributed, $\Eq = O(\log N / \log\log N)$.
Furthermore, there are distributions (in which $X$ and $Y$ are independent and uniformly distributed over sets
of numbers having only small and distinct prime factors) for which we have
$\Eq = \Me(\log N / \log\log N)$.
\vfill\eject
}

\heading{1.  Introduction}

Let $X$ and $Y$ be random integers.
We regard a multiplier as 
a deterministic channel whose input is the pair $(X,Y)$ and whose output is the product $\XY$.
The information lost in multiplication is, according to Shannon [S3], the equivocation $H(X,Y\mid \XY)$.
From the definition of conditional entropy, we have
$$\eqalignno{
\Eq &= H(X,Y,\XY) - H(\XY) \cr
& = H(X,Y) - H(\XY),  &(1.1)\cr
}$$
where we have used the fact that the channel is deterministic ($X\cdot Y$ is determined by $X$ and $Y$,
so that $ H(X,Y,X\cdot Y) = H(X,Y)$).

We first consider the case in which $X$ and $Y$ are independent and uniformly distributed on the set $\{1,\ldots, N\}$,
so that $H(X,Y) = 2\log_2 N$.
We shall show in Section 2 that in this case we have
$$\Eq = \Th(\log\log N). \eqno(1.2)$$

If $X$ and $Y$ have arbitrary (that is, not necessarily independent or uniform)  distributions on $\{1,\ldots, N\}$,
then it is of course possible that $\Eq = 0$.
We may then ask how close $H(X,Y)$ can come to its maximum $2\log_2 N$, while still achieving
$\Eq = 0$. 
We shall show in Section 3 that $\Eq = 0$ implies that
$$2\log_2 N - H(X,Y) = \Me(\log\log N). \eqno(1.3)$$
Furthermore, by taking $X$ and $Y$ to be independent, with distributions 
concentrated on disjoints sets of primes, 
it is possible
to achieve $\Eq = 0$ with $\log_2 N - H(X)$ and $\log_2 N - H(Y)$ each $O(\log\log N)$, so that (1.3) is the
best possible bound.

We shall also consider the distributions of $X$ and $Y$ that maximize the information loss.
We shall show in Section 4 that for any distributions of $X$ and $Y$ on $\{1,\ldots, N\}$ we have
$$\Eq = O(\log N / \log\log N). \eqno(1.4)$$
Furthermore, by taking $X$ and $Y$ to be independent, with distributions
concentrated  on integers having only small and distinct prime
factors, we can achieve
$$\Eq = \Me(\log N / \log\log N),$$
so that (1.4) is the best possible bound.
\vfill\eject

Results concerning information flow through a multiplier have been used by Abelson and Andreae [A] and by
Brent and Kung [B] to obtain
lower bounds involving the area and time required for multiplication.
Furthermore, the results in Section 4 give a lower bound to the number of ancillary lines required by
a reversible multiplier (see Fredkin and Toffoli [F] for a discussion of reversible computation).
This lower bound is achievable if multiplication is performed by a single gate; it is an open question whether
it can be achieved if the multiplier is implemented using standard reversible gates, such as those
proposed by Fredkin and Toffoli.

The proofs in this paper draw upon a variety of results from number theory.
Many of these in turn rely on the prime-number theorem (first proved by Hadamard [H1] and independently
by de la Vall\'{e}e Poussin [V]) and its extension to primes in arithmetic progressions
(first proved by de la Vall\'{e}e Poussin [V]).
While these deep theorems now have elementary proofs (due to Selberg [S1, S2] and Erd\H{o}s [E1]),
none of our results actually depend on theorems of this depth, and thus we shall take care to point out
the simplest results that support our proofs.
\sk

\heading{2. The Uniform Distribution}

Our goal in this section is to establish (1.2).
For $X$ and $Y$ independent with the uniform distribution, we have
$$H(X,Y) = 2\log_2 N.$$
Thus from (1.1) we have
$$\Eq = 2\log_2 N - H(\XY). \eqno(2.1)$$

Define $m(N)$ by
$$m(N) = \#\{x\cdot y : 1\le x\le N, 1\le y\le N\}.$$
We have 
$$H(\XY) \le \log_2 m(N).$$
Thus the bound
$$\Eq = \Me(\log\log N) \eqno(2.2)$$
is a consequence of (2.1) and the following result.

\label{Proposition 2.1:}
For any $\ep > 0$, we have
$$m(N) \le {N^2 \over (\log N)^{\al - \ep}} \eqno(2.3)$$
for all sufficiently large $N$, where 
$\al = 1 - \log_2 (e\ln 2) = 0.08607\ldots\,$.
\vfill\eject

This result is due to Erd\H{o}s [E2], who also proved the matching bound
$$m(N) \ge {N^2 \over (\log N)^{\al + \ep}}.$$
For completeness, we shall give a simple proof of this proposition.

\label{Proof of Proposition 2.1:}
Let $f(n)$ denote the number of distinct prime factors in the integer $n\ge 1$.
Let $\ta_k(x)$ denote the number of integers $n$ in the interval $1\le n\le x$ 
such that  $f(n) = k$.
Hardy and Ramanujan [H2] (Lemma B) have shown that there are absolute constants $L$ and $D$ such that
$$\ta_k(x) \le {L x \over \ln x}\,{(\ln\ln x + D)^{k-1} \over (k-1)!} \eqno(2.4)$$
for all $k\ge 1$ and $x\ge 2$.
Apart from an elementary precursor
$$\ta_1(x) = O\left({x\over \log x}\right)$$
to the prime-number theorem due to Chebyshev [C], their result relies only on the elementary estimates
$$\sum_{p\le x}{\ln p\over p} = \ln x + O(1)$$
and
$$\sum_{p\le x} {1\over p} = O(\log\log x)$$
(in which the sums are over primes $p$) due to  Mertens [M].
We observe that (2.4) implies
$$\ta_k(x) \le {M x \over \ln x}\,{(\ln\ln x)^{k-1} \over (k-1)!} \eqno(2.5)$$
for all $x\ge 2$ and $1\le k\le 2\log_2\ln x$, where $M = L\exp(2D\log_2 e)$.

Fix $0 < \de < 1/6$.
Define $m_1(N)$, $m_2(N)$ and $m_3(N)$ by
$$\leqalignno{
m_1(N) &= \#\big\{(x, y) : 1\le x\le N, 1\le y\le N \hbox{\ and\ } f(x) + f(y) \le (1+2\de)\log_2 \ln N\big\}, \cr
m_2(N) &= \#\big\{z : 1\le z\le n^2 \hbox{\ and\ } f(z)\ge (1+\de)\log_2 \ln N\big\}  \cr
&&\hbox{\ and\ }\cr
m_3(N) &= \#\big\{z : 1\le z\le N^2 \hbox{\ and $w^2 \mid z$ for some $w$ with $f(w)\ge \de\log_2 \ln N$}\big\}. \cr 
}$$
Then we have
$$m(N)\le m_1(N) + m_2(N) + m_3(N).$$
For if $z = x\cdot y$ is not counted by $m_1(N)$, then we have $f(x) + f(y) > (1+2\de)\log_2 \ln N$.
If in addition $z$ is not counted by $m_2(N)$, then we have $f(z) < (1+\de)\log_2 \ln N$, and thus
$$f\(\gcd(x, y)\) > \de\log_2\ln N,$$
where $\gcd(x,y)$ denotes the greatest common divisor of $x$ and $y$.
If we now let $w$ be the product of the distinct primes dividing $\gcd(x,y)$, then we have $w^2\mid z$
and $f(w)\ge \de\log_2 \ln N$, so that $z = x\cdot y$ is counted by $m_3(N)$.
Thus it will suffice to show that $m_1(N)$, $m_2(N)$ and $m_3(N)$ each satisfy a bound of the form of that in (2.3). 

For $m_1(N)$ we have
$$\eqalignno{
m_1(N) &\le \sum_{i+j\le (1+2\de)\log_2 \ln N} \ta_i(N)\,\ta_j(N) \cr
&\le {M^2\,N^2\over (\ln N)^2}\sum_{i+j=k\le (1+2\de)\log_2 \ln N} {(\ln\ln N)^{k-2} \over (i-1)!\,(j-1)!} \cr
&\le {M^2\,N^2\over (\ln N)^2}\sum_{i+j=k\le (1+2\de)\log_2 \ln N} {k-2\choose i-1}{(\ln\ln N)^{k-2} \over (k-2)!} \cr
&\le {M^2\,N^2\over (\ln N)^2}\sum_{i+j=k\le (1+2\de)\log_2 \ln N} {(2\ln\ln N)^{k-2} \over (k-2)!} \cr
&\le {M^2\,N^2\over (\ln N)^2}\sum_{i+j=k\le (1+2\de)\log_2 \ln N} {\left(2e\ln\ln N \over k-2\right)^{k-2}}, &(2.6)\cr
}$$
where we have used the definition of $m_1(N)$, the bound (2.5), the identity $a!/b!\,(a-b)! = {a\choose b}$,
the inequality ${a\choose b}\le 2^a$ and the inequality $a!\ge a^a/e^a$.

The summand in (2.6) increases with $k$ for $k-2\le 2\ln\ln N$, and decreases thereafter.
Since $k-2 < (1+2\de)\log_2 \ln N \le 2e\ln\ln N$,
the largest terms of the sum are those with the largest $k$.
There are at most $\((1+2\de)\log_2 \ln N\)^2 \le (2\ln\ln N)^2$ terms in all,
and each term is at most
$$(2e\ln 2)^{(1+2\de)\log_2\ln N} = (\ln N)^{(1+2\de)(2-\al)}.$$
Thus we obtain the bound
$$m_1(N) \le {M^2 \, (2\ln\ln N)^2 \, (\ln N)^{2\de(2-\al)} \, N^2 \over (\ln N)^\al},$$
which is of the form desired, since if $2\de(2-\al) < \ep$,
the factors 
$$M^2 \, (2\ln\ln N)^2 \, (\ln N)^{2\de(2-\al)}$$
in the numerator
can be absorbed by the factor $(\ln N)^\ep$ in the denominator of (2.3).
\vfill\eject

For $m_2(N)$ we have
$$\eqalignno{
m_2(N) &\le \sum_{k\ge (1+\de)\log_2 \ln N} \ta_k(N^2) \cr
&\le {M\, N^2 \over \ln N} \sum_{k\ge (1+\de)\log_2 \ln N} {(\ln\ln N)^{k-1} \over (k-1)!} \cr
&\le {M\, N^2 \over \ln N} \sum_{k\ge (1+\de)\log_2 \ln N} \left({e\ln\ln N \over k-1}\right)^{k-1}, &(2.7)\cr
}$$
where we have used the definition of $m_2(N)$, the bound (2.5) 
the inequality $a!\ge a^a/e^a$.

The summand in (2.7) increases with $k$ for $k-1\le \ln\ln N$, and decreases thereafter.
Since $k-1 \ge (1+\de)\log_2 \ln N  - 1\ge \ln\ln N$,
the largest terms of the sum are those with the smallest $k$.
There are at most $2e\log_2 \ln N$ terms with $k-1 < 2e\log_2 \ln N$,
and each such term is at most
$$(e\ln 2)^{(1+\de)\log_2\ln N} = (\ln N)^{(1+\de)(1-\al)}. \eqno(2.8)$$
Furthermore, all the terms with $k-1\ge 2e\log_2 \ln N$ are bounded by the terms of a geometric progression
with ratio $1/2$, and thus their sum is bounded by (2.8).
Thus we obtain the bound
$$m_2(N) \le {M \, (1 + 2e\log_2\ln N) \, (\ln N)^{\de(1-\al)} \, N^2 \over (\ln N)^\al},$$
which is of the form desired, since if $\de(1-\al) < \ep$,
the factors 
$$M \, (1 + 2e\log_2\ln N) \, (\ln N)^{\de(1-\al)}$$ 
in the numerator
can be absorbed by the factor $(\ln N)^\ep$ in the denominator of (2.3).

Finally, for $m_3(N)$ we have
$$\eqalignno{
m_3(N) &\le \sum_{f(w)\ge \de\log_2 \ln N} {N^2\over w^2} \cr
&\le  \sum_{w\ge w_0} {N^2 \over w^2} \cr
&\le {N^2 \over w_0}, &(2.9) \cr
}$$
where $w_0$ denotes the smallest integer $w$ such that $f(w)\ge \de\log_2 \ln N$.
Clearly $w_0 = p_1 \cdots p_k$ is the product of the first $k = \lceil \de\log_2 \ln N\rceil$ primes.
If $N$ is sufficiently large that there are fewer than $k/2$ primes that are less than $2^{2/\de}$,
then $w_0$ contains at least $k/2$ prime factors that are each at least $2^{2/\de}$, and thus
$w_0 \ge \ln N$.
The bound (2.9) is therefore also of the desired form.
This completes the proof of the proposition.
\QED
\vfill\eject

Next we turn to establishing the upper bound
$$\Eq = O(\log\log N). \eqno(2.10)$$
To do this we use the formula
$$\eqalignno{
\Eq &= \sum_{1\le x\le N}  \sum_{1\le y\le N} \Pr[X=x,Y=y] \, H(X,Y\mid \XY = \xy). &(2.11)\cr
}$$ 
Using the bound
$$\eqalign{
H(X,Y\mid \XY = \xy)
&\le \log_2 \#\big\{ (v,w) : 1\le v\le N, 1\le w\le N\hbox{\ and\ } v\cdot w = \xy\big\} \cr
&\le \log_2 d(\xy) \cr
}$$
(where $d(n)$ denotes the number of divisors of the integer $n$), we obtain
$$\Eq \le \sum_{1\le x\le N}  \sum_{1\le y\le N} \Pr[X=x,Y=y] \, \log_2 d(\xy). \eqno(2.12)$$

For $X$ and $Y$ independent with
the uniform distribution, (2.12) becomes
$$\Eq \le {1\over N^2} \sum_{1\le x\le N}  \sum_{1\le y\le N}  \log_2 d(\xy). \eqno(2.13)$$
Since $\log_2 a$ is a concave function of $a$, the average of the logarithm in (2.13) is at most the 
logarithm of the average, and we obtain
$$\Eq \le \log_2 \left({1\over N^2} \sum_{1\le x\le N}  \sum_{1\le y\le N}   d(\xy)\right).$$
Since $d(\xy) \le d(x)\cdot d(y)$, we obtain
$$\eqalignno{
\Eq
&\le \log_2 \left({1\over N^2} \sum_{1\le x\le N} \sum_{1\le y\le N} d(x)\cdot d(y)\right) \cr
&=  2\log_2 \left({1\over N} \sum_{1\le n\le N} d(n)\right). &(2.14)\cr
}$$
We now use the asymptotic formula
$$\sum_{1\le n\le N} d(n) = N\ln N + O(N)$$
due to Dirichlet [D2] (which is established simply by estimating the number of lattice points in the 
region bounded by the $x$-axis, the $y$-axis and the hyperbola $\xy = N$).
Substituting this result in (2.14) completes the proof of (2.10), which together with (2.2) establishes (1.2).
\vfill\eject

\heading{3.  Multiplication without Loss of Information}

Our goal in this section is to determine the maximum entropy that $X$ and $Y$ can have when
$\Eq = 0$.
Let 
$$\calW = \{(x,y) : \Pr[X=x,Y=y] > 0\}$$
denote the support of the distribution
of $(X,Y)$, and let
$$\calM = \{x\cdot y : 1\le x\le N, 1\le y\le N\}$$
be the range of the multiplication
map $\mu : \{1, \ldots, N\}\times \{1, \ldots, N\}\to \{1, \ldots, N^2\}$ defined by
$\mu(x,y) = \xy$.
Then $\Eq = 0$ implies that $\mu$ restricted to $\calW$ is injective, so that
$\#(\calW) \le \#\calM = m(N)$ and $H(X,Y)\le \log_2 m(N)$.
Proposition 2.1 thus shows that $\Eq = 0$ implies (1.3).

To show that this result is the best possible, we let $X$ and $Y$ be 
independent and uniformly distributed 
over $\calX$ and $\calY$, respectively, where $\calX$ and $\calY$
are the sets of primes that are
at most $N$ and congruent to $1$ and $3$, respectively, modulo $4$. 
To show that $\log_2 N - H(X)$ and $\log_2 N - H(Y)$ are each $O(\log \log N)$, it will suffice
to show that $\#\calX = \pi_{1,4}(N)$ and $\#\calY = \pi_{3,4}(N)$ are each $\Me(N/\log N)$.
This of course follows from the extention of the prime-number theorem to arithmetic progressions,
but we can obtain what we need from the following simple result due to Shapiro [S4]
(which is an elementary quantitative version of the theorem of Dirichlet [D1] on primes in 
arithmetic progressions).
Let $a$ and $b$ be fixed with $\gcd(a,b)=1$.
Then
$$\sum_{\textstyle p\le x \atop \textstyle p\equiv a \,({\rm mod}\, b)} {\ln p\over p} 
= {\ln x\over \ph(b)} + O(1), \eqno(3.1)$$
where $\ph(b)$ denotes Euler's totient function: the number of $a$ in the range $0 < a < b$ such 
$\gcd(a,b)=1$.
To show that (3.1) implies
$$\pi_{a,b}(x) = \Me\left({x\over \log x}\right),  \eqno(3.2)$$
we observe that (3.1) implies that
$$\sum_{\textstyle x/A < p\le x \atop \textstyle p\equiv a \,({\rm mod}\, b)} {\ln p\over p} 
\ge {\ln A\over \ph(b)} - 2B \eqno(3.3)$$
for all $A > 1$, where $B$ is a bound on the magnitude of the $O(1)$ term in (3.1).
Choosing $A$ sufficiently large that the right-hand side of (3.3) is strictly positive and
observing that each term in the sum is at most $(A\ln x) / x$ establishes that there must be 
$\Me(x/\log x)$ terms, and thus yields (3.2).
\vfill\eject

\heading{4.  The Maximum Loss of Information}

Our goal in this section is to determine the maximum possible loss of information in multipication.
Our starting point is the formula (2.12).
Since the average is at most the maximum, we have
$$\Eq \le \max_{1\le x\le N}  \max_{1\le y\le N} \log_2 d(\xy),$$
and since $\log_2 a$ is an increasing function of $a$, we obtain
$$\Eq \le \log_2 \left(\max_{1\le x\le N}  \max_{1\le y\le N}  d(\xy)\right).$$
Using the fact that $d(\xy)\le d(x)\cdot d(y)$ as before,
we obtain
$$\Eq \le 2\log_2 \left(\max_{1\le n\le N} d(n)\right). \eqno(4.1)$$
Wigert [W] was the first to show that
$$\log_2 \left(\max_{1\le n\le N} d(n)\right) \sim {\ln N\over \ln\ln N}, \eqno(4.2)$$
using the prime-number theorem.
But Ramanujan [R] has shown that 
an estimate even more precise than
(4.2) can be obtained
using only the crude bounds
$$\pi(x) = \Th\left({x\over\log x}\right) \eqno(4.3)$$
for the number $\pi(x)$ of primes not exceeding $x$ obtained by Chebyshev [C].
Substituting (4.2) into (4.1) yields (1.4).

To show that this result is the best
possible, we let $X$ and $Y$ be independent and uniformly distributed on the set $\calV$  of the $2^k$ divisors
of the product $v_k = p_1 \cdots p_k$ of the first $k$ primes, where $k$ is the largest integer such that
$$v_k \le N. \eqno(4.4)$$
If we define $\th(x)$ by
$$\th(x) = \sum_{p\le x} \ln p$$
(in which the sum is over primes $p$), then 
$$v_k = \exp \th(P_k),$$
so that (4.4) is equivalent to 
$$\th(p_k) \le \ln N.$$
The bounds
$$\th(x) = \Th(x)$$
are equivalent to the bounds (4.3) established by Chebyshev [C].
This implies that 
$$p_k = \Th(\log N),$$
so that (again using (4.3))
$$k = \Th\left({\log N\over \log\log N}\right). \eqno(4.5)$$

From (2.11), we have
$$\Eq = {1\over 2^{2k}}\sum_{x\in\calV} \sum_{y\in\calV} H(X,Y\mid \XY=\xy). \eqno(4.6)$$
For $x, y\in\calV$, let $u(x,y)$ denote the number of primes among $p_1, \ldots, p_k$ that divide one, 
but not both, of $x$ and $y$.
(This number is also the number of primes that divide the square-free part of $\xy$,
and thus it depends only on $\xy$.)
The random variable $(X,Y)$, conditioned on $\XY=\xy$, is uniformly distributed over the $2^{u(x,y)}$
pairs in the set
$$\calU =\{(v,w)\in\calV\times \calV : v\cdot w = \xy\},$$
so that
$$\eqalign{
H(X,Y\mid \XY = \xy) &= \log_2 \#\calU \cr
&= u(x,y). \cr}$$
Thus (4.6) yields
$$\Eq = {1\over 2^{2k}}\sum_{x\in\calV} \sum_{y\in\calV} u(x,y). \eqno(4.7)$$

Since $X$ and $Y$ are each uniformly distributed on the $2^k$ divisors of $v_k$, 
the divisibility of each of $X$ and $Y$ by each of the primes $p_1, \ldots, p_k$ is probabilistically equivalent to
the occurrences of heads among $2k$ independent flips of an unbiased coin. 
In particular, each of the primes $p_1, \ldots, p_k$ divides one, but not both, of $X$ and $Y$ with
probability $1/2$.
Thus the right-hand side of (4.7) is equal to $k/2$, and (4.5) yields
$$\eqalign{
\Eq &= k/2 \cr
&= \Me\left({\log N\over \log\log N}\right). \cr
}$$
This estimate shows that the result (1.4) is the best possible.
\sk

\heading{5.  References}

\ref A; H. Abelson and P. Andreae;
``Information Transfer and Area-TimeTrade-Offs for VLSI Multiplication'';
Comm.\ ACM; 23 (1980) 20--23.

\ref B; R. P. Brent and H. T. Kung;
``The Area-Time Complexity of Binary Multiplication'';
Journal of the Association for Computing Machinery; 28 (1981) 521--534;
Corrigendum: 29 (1982) 904.

\ref C; P. L. Chebyshev ($=$ Tchebichef);
``M\'{e}moire sur les nombres premiers'';
Journal de math\'{e}matiques pures et appliqu\'{e}es (1); 17 (1852) 366--390.

\ref D1; P. G. L. Dirichlet ($=$ Lejeune-Dirichlet);
``Sur l'usage des s\'{e}ries infinies  dans la th\'{e}orie des nombres'';
Journal f\"{u}r die reine und angewandte Mathematik''; 18 (1838) 259--274.

\ref D2; P. G. L. Dirichlet ($=$ Lejeune-Dirichlet);
``Valeurs moyennes dans la th\'{e}orie des nombres'';
Journal de math\'{e}matiques pures et appliqu\'{e}es (2); 1 (1956) 353--370.

\ref E1; P. Erd\H{o}s;
``On a New Method in Elementary Number Theory Which Leads to an  Elementary Proof of the Prime Number Theorem'';
Proceedings of the National Academy of Science of the USA; 35 (1949) 374--384.

\ref E2; P. Erd\H{o}s;
``Ob Odnom Asimptoticheskom Neravenstve v Teorii Chisel'';
Vestnik Leningradskogo Universiteta; 13 (1960) 41--49.

\ref F; E. Fredkin and T. Toffoli;
``Conservative Logic'';
International Journal of Theoretical Physics; 21 (1982) 41--55.

\ref H1; J. Hadamard;
``Sur la distribution des z\'{e}ros de la fonction $\ze(s)$ et ses cons\'{e}quences arithm\'{e}tiques'';
Bulletin de la Soci\'{e}t\'{e} Math\'{e}matiques de France''; 24 (1896) 199--220.

\ref H2; G. H. Hardy and S. Ramanujan;
``The Normal Number of Prime Factors of a Number $n$'';
Quarterly Journal of Mathematics; 48 (1917) 76--92.

\ref M; F. Mertens;
``Ein Beitrag zur analytischen Zahlentheorie'';
Journal f\"{u}r die reine und angewandte Mathematik''; 78 (1874) 46--62.

\ref R; S. Ramanujan;
``Highly Composite Numbers'';
Proceedings of the London Mathematical Society (2); 14 (1915) 347--409.

\ref S1; A. Selberg;
``An Elementary Proof of the Prime-Number Theorem'';
Annals of Mathematics; 50 (1949) 305--313.

\ref S2; A. Selberg;
``An Elementary Proof of the Prime-Number Theorem for Arithmetic Progressions'';
Canadian Journal of Mathematics; 2 (1950) 66--78.

\ref S3; C. E. Shannon;
``A Mathematical Theory of Communication'';
Bell System Technical Journal; 27 (1948) 379--423, 623--655.

\ref S4; H. N. Shapiro;
``On Primes in Arithmetic Progression (II)'';
Annals of Mathematics; 52 (1950) 231--243.

\ref V; Ch.\ de la Vall\'{e}e Poussin;
``Recherches analytiques sur la th\'{e}orie des nombre premiers'';
Annales de la Soci\'{e}t\'{e} Scientifique de Bruxelles; 20 (1896) 183--256, 281--397.

\ref W; S. Wigert;
``Sur l'ordre de grandeur du nombre des diviseurs d'un entier'';
Arkiv f\"{o}r Matematik, Astronomi och Fysik; 3, 18 (1907) 1--9.

\bye